# On the local Langlands correspondence mod $\ell$

Chandrashekhar Khare

## 1 Introduction

We prove in this paper that supercuspidal representations of $GL_n(F)$, with $F$ a finite extension of $\mathbf{Q}_p$, that are congruent mod $\ell$ ($\ell \neq p$) have Langlands parameters that, up to semisimplification, are congruent mod $\ell$.

To state the main theorem precisely, we fix embeddings of $\overline{\mathbf{Q}}$ into $\overline{\mathbf{Q}_p}$ and $\overline{\mathbf{Q}_\ell}$ and denote the corresponding places of $\overline{\mathbf{Q}}$ by $\wp$ and $\ell$ respectively. We also fix an isomorphim $\overline{\mathbf{Q}_\ell} \cong \mathbf{C}$, and after this may view admissible representations over $\mathbf{C}$ and $\overline{\mathbf{Q}_\ell}$ equivalently, as admissible representations do not notice the topology of the coefficient field.

The theorem proved in this paper is:

**Theorem 1** *Assume that $\ell \neq p$. If $\pi_\wp$ and $\pi'_\wp$ are supercuspidal representations of $GL_n(K)$ over $\overline{\mathbf{Q}_\ell}$, with identical finite order central character $\omega$, and such that the mod $\ell$ reductions $r_\ell(\pi_\wp)$ and $r_\ell(\pi'_\wp)$ are isomorphic, then the semisimplification of the mod $\ell$-reduction of the corresponding Galois representations (i.e., Langlands parameters) $\sigma_\wp$ and $\sigma'_\wp$ are isomorphic.*

**Remarks.**

1. Vigneras has proven that under the assumption of the theorem $\pi_\wp$ (resp. $\pi'_\wp$) has an integral model and its reduction mod $\ell$, denoted by $r_\ell(\pi_\wp)$ (resp., $r_\ell(\pi'_\wp)$), is well-defined and in fact irreducible; this is because of the assumption of supercuspidality.

2. In the preprint [V] the theorem above was proven for the case of $\ell > n$ and $\ell \neq p$, for $\ell > n$. This result was a key step in loc. cit. to prove that the classical local Langlands correspondence induces a map between $\ell$-supercuspidal representations of $GL_n(K)$ (over $\overline{\mathbf{Q}_\ell}$), and $F$-semsimple representations of $WD_K$, the Weil-Deligne group of $K$, that are irreducible on reduction mod $\ell$. These are called $\ell$-irreducible parameters. It is a general fact (see [V1]) that $r_\ell(\pi_\wp)$ is irreducible and *cuspidal* when $\pi_\wp$ is supercuspidal: but it need not be *supercuspidal*. The methods of the present paper (that are global) and those of [V] (using local harmonic analysis) are completely different.

3. Both the results and methods of the present paper (in particular the proof of Lemma 1 below) have been used by Vigneras in a recent preprint [V2] to give a different proof using global methods of the local Langlands correspondence for $GL_n(F)$ mod $\ell$ for $F$ any non-arhimedean local field of residue characteristic $p$ with $\ell \neq p$. The local methods used in [V] did not work when $F$ was of characteristic $p$, and nor when $\ell < n$.

4. We note the general fact that the reduction mod $\ell$ is well-defined (i.e., independent of the choice of lattice) upto semisimplification. For admissible representations of $GL_n(K)$ this is in [V1], while on the Galois side this is Brauer-Nesbitt.

5. In the proof below we need more than just the abstract statement of the local Langlands correspondence, and make use of the results of Harris and Taylor, rather than Henniart, cf. [H]. Namely we need the compatibility proven proved in [HT] in many cases between the local Langlands correspondence, and the restriction at $\wp$ of the $\ell$-adic representations attached to self-dual automorphic representations on $GL_n$ over totally real fields; these proven cases suffice to prove Theorem 1. In the preprint [V2] the role of the results of Harris and Taylor is filled by the recent work of Lafforgue [L] on the Langlands conjecture in the function field case.

## 1.1 Sketch of proof

The idea of the proof is to embed $\pi_\wp$ and $\pi'_\wp$ global automorphic representations that are cohomological and which are congruent mod $\ell$, and for these to use the information at places at which the representations are unramified principal series to conclude. By two global representations being congruent mod $\ell$, here and throughout this paper, we just mean a congruence of the Satake parameters at all places where both the representations are unramified modulo the place of $\overline{\mathbf{Q}}$ fixed by the embedding $\overline{\mathbf{Q}} \to \overline{\mathbf{Q}_\ell}$.

Fix a totally real number field $F^+$ that localises to $K$ at a place of $F^+$ that is induced by $\wp$, and that we denote by abuse of notation by the same symbol $\wp$. In a little more detail, we seek automorphic representations $\Pi$ and $\Pi'$ in the cohomology of Kottwitz varieties attached to unitary groups over $F^+$, arising from division algebras with an involution of the second kind over a CM extension $F$ of $F^+$, that are congruent mod $\ell$, and such that $\Pi_\wp = \pi_\wp$ and $\Pi'_\wp = \pi'_\wp$. Work of Kottwitz, Clozel and Taylor (cf. [Ha]) attaches a $\ell$-adic Galois representation $\Sigma$ and $\Sigma'$ of $Gal(\overline{F}/F)$ to $\Pi$ and $\Pi'$, that is related in the usual way that is recalled below. Then we deduce that the semsimplifications of the mod $\ell$-reductions of $\Sigma|_{D_\wp}$ and $\Sigma'|_{D_\wp}$ are isomorphic as $\Pi$ and $\Pi'$ are congruent mod $\ell$. Harris and Taylor prove that the association $\pi_\wp \to \Sigma|_{D_\wp}$ (resp., $\pi'_\wp \to \Sigma'|_{D_\wp}$) *is* the local Langlands correspondence upto taking duals and twisting (the prescription is recalled in the proof below; see also page 26 of [Ca]). Thus on reducing mod $\ell$ we conclude that $\sigma_\wp$ and $\sigma'_\wp$ have isomorphic mod $\ell$ reductions upto semisimplification.

The method to find a cohomological $\Pi$ with local component $\pi_\wp$ is by now a standard application of the trace formula (see [Cl] and [Ha]). To find a congruent $\Pi'$ with local component $\pi'_\wp$ is in principle an application of a suitable version of Carayol's lemma (cf. Lemme 2 of [Ca1] and Section 2 of [DT]). But as we do not know torsion-free properties of cohomology of Kottwitz varieties we first work with unitary groups that are compact mod centre at infinity. Here the analog of Carayol's lemma is easy (this is Lemma 1 below). Then we use Lemma 5 of [Ha], that in turn is based on the pseudo-stabilisation of the trace formula due to Kottwitz and Clozel (see also [Cl1]), to transfer $\Pi$ and $\Pi'$ to a unitary group that arises from a division algebra with involution of the second kind, but that is now of the type considered in [HT]: in particular it has type $(n-1, 1)$ at one infinite place and $(n, 0)$ at all other places.

## 2 The proof

We fix a number field $F$ that is the composite of a totally real field $F^+$ and an imaginary quadratic extension $\mathbf{Q}(\sqrt{-D})/\mathbf{Q}$ in which the prime $p$ splits.

The relationship between $F$ and $K$ is that the localisation of $F^+$ at $\wp$ is $K$. Let $\wp_1$ and $\wp_2$ be the places of $F$ above $\wp$; thus the localisation of $F$ at $\wp_i$ is again $K$.

Let $B$ be a division algebra of degree $n$ over $F$, with an involution $*$ that induces the non-trivial automorphism of $F$ over $F^+$, and such that it is either a division algebra or split at any place of $F$. Let $\beta \in B$ be a $*$-antisymmetric element, and conjugating $*$ by $\beta$ one obtains another involution that we denote by $\iota$. We assume that $B$ splits at $\wp_1$ and $\wp_2$, and is ramified at sufficiently many other finite places: this will derive meaning from the considerations below.

We associate a fake unitary group $G$ over $F^+$ to the pair $(B, \iota)$ such that

$$G(R) = \{g \in B^{op} \otimes R; b\iota(b) \in R^*\}$$

for any $F^+$-algebra $R$. We assume that $G(F^+ \otimes \mathbf{R})$ is of type $(n, 0)$ at all the infinite places. A pair $(B, \iota)$ that gives rise to such a $G$ exists by the considerations in [Cl1] related to the Hasse principle as we are willing to allow $B$ (resp., $G$) to be ramified at a sufficiently large set of places of $F$ (resp., $F^+$).

Using the same arguments, and the fact that $G$ is sufficiently ramified, we can deduce there exists a division algebra $B'$ with an involution $\iota'$ of the second kind over $F$ with the following properties:

1. It is either split or a division algebra at all places.

2. It is split at $\wp_1$ and $\wp_2$.

3. The unitary group $G'$ associated to the pair $(B', \iota')$ has type $(n-1, 1)$ at one infinite place and type $(n, 0)$ at the other infinite places of $F^+$.

4. The finite places of $F^+$ at which $G'$ is not split are contained in the finite places of $F^+$ at which $G$ is not split.

As $\pi_\wp$ has finite order central character $\omega$, then by the result recorded as Lemma 1 of [Ha], and due to Clozel [Cl], we may find an automorphic representation $\Pi$ of $G(\mathbf{A}_{F^+})$ that has local component $\pi_\wp$ at $\wp$, and is the trivial representation at the infinite places.

We use the notation of the statement of Theorem 1. The following lemma is the key technical result of the paper:

**Lemma 1** *There is an automorphic representation $\Pi'$ of $G(\mathbf{A}_{F^+})$, with local component $\pi'_\wp$ at $\wp$, that is congruent to $\Pi$ mod $\ell$, and that is the trivial representation at infinity.*

**Proof.** We assume that $n > 1$ as the case $n = 1$ is straightforward. We begin by recalling the theory of types of Bushnell-Kutzko for admissible supercuspidal representations of $GL_n(K)$ following the exposition in [V] Section 2.4 and in [V1] III.4.26.

It follows from the work of Bushnell-Kutzko that every supercuspidal representation $V$ of $GL(n)$ is obtained from an irreducible representation $\Lambda$ of a maximal compact-mod-centre subgroup $J$ of $GL(n)$ by (compact) induction, such that the isomorphism class of the pair $(J, \Lambda)$ is uniquely determined by $V$ up to conjugation action on this pair by $GL(n)$. Vigneras has proved (loc. cit.) that if two supercuspidal representations of $GL(n)$ have the same reduction mod $\ell$, which by another theorem of hers is irreducible, then we can choose the corresponding $J$'s to be the same and further the associated $\Lambda$'s will have the same reduction mod $\ell$. We denote by $(J, \Lambda)$ and $(\overline{J}, \Lambda')$ the data associated to $\pi_\wp$ and $\pi'_\wp$: as we are assuming that the latter are congruent mod $\ell$ we can take $J = \overline{J}$ and the (irreducible) reductions $r_\ell(\Lambda)$ and $r_\ell(\Lambda')$ are isomorphic. Note that these are of dimension $> 1$.

We now define the space of "modular forms" that are useful for us. This is done in a way analogous to Section 2 of [DT]. We fix a sufficiently deep open compact-mod-centre subgroup $U$ of $G(\mathbf{A}^f_{F^+})$ that at $\wp$ is just the $J$ above. By sufficiently deep we mean that it satisfies (i) and (ii) below:

(i) If $G(\mathbf{A}_{F^+,f}) = \cup_i G(F^+) g_i U$, the disjoint union of finitely many double cosets $G(F^+) g_i U$, then $G(F^+) \cap g_i U g_i^{-1}$ is central for all $i$. This follows for $U$ deep enough as this intersection is finite mod the centre (in $G(F^+)$) which in turn follows from

(a) $G$ is compact mod centre at infinity

(b) For any sufficiently large prime $q$ of $K$ the subgroup of matrices of $GL_n(\mathcal{O}_q)$, with $\mathcal{O}_q$ the ring of integers of $K_q$, congruent to 1 mod $q$ is torsion-free.

(ii) We also require that $\Pi^{U^{(\wp)}} \neq 0$, where by $U^{(\wp)}$ we mean $U$ deprived of its $\wp$-component.

The space of modular forms relevant to us is the $\overline{\mathbf{Z}}_\ell$-module

$$\mathcal{X}_{\pi_\wp} := Maps_U(\Lambda \to G(F) \backslash G(\mathbf{A}) / G_\infty)$$

and the analogously defined $\mathcal{X}_{\pi'_\wp}$. We explain the notation. By $\Lambda$ we mean a $\overline{\mathbf{Z}_\ell}$- model of $\Lambda$ on which $U$ acts through its component at $\wp$ which is $J$. The subscript $U$ denotes $U$-equivariant maps. Note that we do not have to mod out by the *trivial* subspace as in [DT] because $\Lambda$ is irreducible and of dimension bigger than 1.

By our choice of $U$ it satisfies property (i) above, and thus we can view $\mathcal{X}_{\pi_\wp}$ as just the sum of copies of $\oplus_i \Lambda^{Z_i}$, where the superscripts mean fixed points under the central subgroups $Z_i = g_i^{-1} G(F^+) g_i \cap U$, indexed by the double coset space $G(F^+) \backslash G(\mathbf{A}_{F^+}) / G_\infty U$.

This space has an action of the Hecke algebra over $\overline{\mathbf{Z}_\ell}$, $\mathcal{H} := \overline{\mathbf{Z}_\ell}[U^S \backslash G(\mathbf{A}_f^S) / U^S]$ where the superscript $S$ denotes a fixed sufficiently large finite set of finite places of $F^+$, that will usually include all those above $p$. We choose $S$ sufficiently large so that $U^S$ is hyperspecial at all places and thus $\mathcal{H}$ is commutative.

For any automorphic representation $\Pi''$ of $G(\mathbf{A}_{F^+})$ that is trivial at infinity, and with the fixed points of $\Pi''_\wp|_J$ under the subgroup of $\hat{J}$ of $J$ that is the kernel of the representation of $J$ associated to $\Lambda$, isomorphic to the $J/\hat{J}$ representation $\Lambda$, with $g \in \Pi''^{U^{(\wp)}}$ ($g \neq 0$) an eigenvector for the Hecke action (away from $S$), there is a non-zero element in $\mathcal{X}_{\pi_\wp}$ with the same eigenvalues for the $\mathcal{H}$-action as $g$. Thus there is a non-zero element in $\mathcal{X}_{\pi_\wp}$ that is an eigenvector for $\mathcal{H}$ and has the same Hecke eigenvalues, or Satake parameters (outside $S$), as a non-zero eigenvector in the space of $\Pi$.

Define
$$\widetilde{\mathcal{X}_{\pi_\wp}} := Maps_U(r_\ell(\Lambda) \to G(F^+) \backslash G(\mathbf{A}) / G_\infty)$$
$$\widetilde{\mathcal{X}_{\pi'_\wp}} := Maps_U(r_\ell(\Lambda') \to G(F^+) \backslash G(\mathbf{A}) / G_\infty)$$
where as before $r_\ell$ denotes the reduction mod $\ell$ that is well-defined upto isomorphism because of Vigneras' result. As we are assuming that $r_\ell(\Lambda) \cong r_\ell(\Lambda')$ (see above), we deduce that $\widetilde{\mathcal{X}_{\pi_\wp}}$ is isomorphic to $\widetilde{\mathcal{X}_{\pi'_\wp}}$ as Hecke modules. Note that by our description above of $\mathcal{X}_{\pi_\wp}$,

$$\mathcal{X}_{\pi_\wp} \otimes_{\overline{\mathbf{Z}_\ell}} \overline{\mathbf{F}_\ell},$$

and

$$\mathcal{X}_{\pi'_\wp} \otimes_{\overline{\mathbf{Z}_\ell}} \overline{\mathbf{F}_\ell},$$

give rise to isomorphic subspaces under the Hecke equivariant isomorphism $\widetilde{\mathcal{X}_{\pi_\wp}} \simeq \widetilde{\mathcal{X}_{\pi'_\wp}}$: for this note that $\pi_\wp$ and $\pi'_\wp$ have the same central character $\omega$.

To have the Hecke equivariant isomorphism

$$\mathcal{X}_{\pi_\wp} \otimes_{\overline{\mathbf{Z}_\ell}} \overline{\mathbf{F}_\ell} \simeq \mathcal{X}_{\pi'_\wp} \otimes_{\overline{\mathbf{Z}_\ell}} \overline{\mathbf{F}_\ell}$$

is the main reason that we have chosen to work with a group $G$ that is compact-mod-centre at infinity.

Thus the maximal ideals (of residue characteristic $\ell$) of $\mathcal{H}$ that are in the support of $\mathcal{X}_{\pi_\wp}$ are the same as those in the support of $\mathcal{X}_{\pi'_\wp}$. This is a consequence of the Deligne-Serre lemma (cf. [DS]).

Now we claim that there is an eigenvector for $\mathcal{H}$ in $\mathcal{X}_{\pi'_\wp}$ that generates an automorphic representation $\Pi'$ congruent to $\Pi$ mod $\ell$ (in the sense described in the introduction above) with local component $\pi'_\wp$. For this note that:

1. $\mathcal{X}_{\pi'_\wp}$ is torsion-free which follows from (the analog of) the explicit description above for $\mathcal{X}_{\pi_\wp}$.

2. If we consider a non-zero eigenvector $f$ in $\mathcal{X}_{\pi'_\wp} \otimes \overline{\mathbf{Q}_\ell}$ for $\mathcal{H}$, then the $G(\mathbf{A}_{F^+,f})$-module generated by $f$ decomposes as a direct sum $\oplus_j \pi_j$ of finitely many irreducible automorphic representations. This follows from our compactness mod centre assumption on $G(F^+ \otimes \mathbf{R})$. Further $\mathcal{X}_{\pi'_\wp} \otimes \overline{\mathbf{Q}_\ell}$ is isomorphic as a Hecke module to

$$[Maps_{U'}(\overline{\mathbf{Q}_\ell} \to G(F) \backslash G(\mathbf{A})/G_\infty)]^{\Lambda'}$$

where $U'$ is the same as $U$ at all places different from $\wp$, and at $\wp$ is the subgroup $J'$ of $J$ that is the kernel of the $J$-representation $\Lambda'$. The space

$$[Maps_{U'}(\overline{\mathbf{Q}_\ell} \to G(F) \backslash G(\mathbf{A})/G_\infty)]$$

has an action of $U/U'$ and by the superscript $\Lambda'$ we mean the $\Lambda'$-isotypical component. From this we see that all the automorphic representations $\pi_j$ have local components $\pi_{j,\wp}$'s at $\wp$, such that their restrictions to $J$ have a subrepresentation that is isomorphic as a $J$ module to $\Lambda'$. This forces all the $\pi_{j,\wp}$'s to be isomorphic to $\pi'_\wp$ (cf. [V1] III.4.27 for instance). To see this note that, as a simple consequence of the existence of types and Frobenius reciprocity, for a supercuspidal representation its "minimal $K$-type" (or in our case more appropriately "minimal $J$-type") occurs in no other non-isomorphic admissible, irreducible representation.

Now as the maximal ideals in the support of $\mathcal{X}_{\pi_\wp}$ and $\mathcal{X}_{\pi'_\wp}$ are the same, we deduce from the above 2 points the existence of a $\Pi'$ as required. Also as noted above $\Lambda'$ is irreducible and of dimension greater than one, and thus $\Pi'$ is forced to be cuspidal of dimension greater than 1. This finishes the proof of Lemma 1.

**Remark.** The existence of $\Pi'$ follows from cohomological arguments, and once $\Pi$ exists one does not need trace formula arguments to create $\Pi'$.

We consider representations of $G'(\mathbf{A}_{F^+})$ that are the functorial lifts of $\Pi$ and $\Pi'$ which exist because of our assumption on the ramification of $G$ and $G'$ by the same arguments as the proof of Lemma 5 of [Ha]. The functorial lift in each case is trivial at all but one infinite place of $F^+$, i.e., at all infinite places where $G'$ is compact, and at the one remaining infinite place is one of the $(U(n-1,1), K_\infty)$ modules $\pi_j$ $(1 \leq j \leq n)$ such that $H^{n-1}(U(n-1,1), K_\infty; \pi_j) \neq 0$, with $K_\infty$ a maximal compact mod centre subgroup of $U(n-1,1)$. We have a lemma that is proved exactly like Lemma 5 of [Ha], and whose proof we omit.

**Lemma 2** *The functorial lifts of $\Pi$ and $\Pi'$ exist in the cuspidal spectrum of $G'$.*

We denote these functorial lifts by the same symbols. By the Matsushima formula $\Pi$ and $\Pi'$ occur in $lim_{\to \Gamma} H^*(\Gamma, \overline{\mathbf{Q}_\ell})$, where $\Gamma$ runs through open compact subgroups of $G'(\mathbf{A}_{F^+}^f)$. Further $* = n-1$ because of our assumption that $\Pi_\wp$ and $\Pi'_\wp$ are supercuspidal (cf. [Cl2]). Consider the following piece of the $\ell$-adic étale cohomology

$$H^{n-1}(\Gamma, \overline{\mathbf{Q}_\ell})[\Pi_f]$$

cut out by the Hecke action for a suitable congruence subgroup $\Gamma$. Then, by work of Kottwitz, Clozel and Taylor, upto semisimplification it is a sum of copies of $\ell$-adic representations $\Sigma(\Pi)$ of $Gal(\overline{F}/F)$ (we called this $\Sigma$ in the introduction and continue to do so below except for a brief while now to emphasise the dependence on $\Pi$) characterised by

$$\sigma_n(\tilde{\Pi}_v) = \Sigma(\Pi)|_{D_v} \otimes |\ |^{(n-1)/2} \qquad (*)$$

where $\sigma_n(\Pi_v)$ is the Langlands parameter of $\Pi_v$ for almost all places $v$ of $F$, and tilde stands for dual. The same holds good for $\Pi'$ and $\Sigma'(:= \Sigma(\Pi'))$.

**Proposition 1** *With $\Pi$ and $\Pi'$, $\Sigma$ and $\Sigma'$ as above, the semisimplifications of the mod $\ell$ reductions of $\Sigma|_{D_\wp}$ and $\Sigma'|_{D_\wp}$ are isomorphic.*

**Remark.** As $\wp$ (regarded as a place of $F^+$) splits into $\wp_1$ and $\wp_2$ in $F/F^+$, we are making identifications of the decomposition groups $D_\wp \cong D_{\wp_i} \cong Gal(\overline{K}/K)$.

**Proof.** Using the relation between $\Pi$ and $\Sigma$ (resp, $\Pi'$ and $\Sigma'$) which is given by the Hecke normalisation (see $(*)$ above), and the fact that $\Pi$ and $\Pi'$ are congruent mod $\ell$, we conclude that the characteristic polynomials of almost all Frobenii in the representations $\Sigma$ and $\Sigma'$ are congruent modulo $\ell$. This makes sense because of the fact that $\Sigma$ and $\Sigma'$ are unramified outside a finite set of places. Note that $\Sigma|_{D_\wp}$ factors through a finite extension of $Gal(\overline{K}/K)$ as $\pi_\wp$ is supercuspidal with finite order central character; this follows from [HT]. The Cebotarev density theorem and the Brauer-Nesbitt theorem (cf. Theorem 30.16 together with Chapter X of [CR]) implies that the semisimplifications of the mod $\ell$ reductions of $\Sigma|_{D_\wp}$ and $\Sigma'|_{D_\wp}$ are isomorphic, proving the proposition.

**Remark.** Thus far we (essentially) did not need [HT] as we have only used information at the good places that was known earlier (see [Cl1]).

It is at this point that we invoke the results of [HT]. As it is precisely in the $\ell$-adic cohomology of these Shimura varieties (i.e., Kottwitz varieties) associated to congruence subgroups $\Gamma$ arising from $G'$ that Harris and Taylor realise the local Langlands correspondence. Namely Harris and Taylor have proven that the correspondence $\pi_\wp \to \Sigma|_{D_\wp}$ *is* the local Langlands correspondence upto the normalisation recalled above in $(*)$, at least when when $\pi_\wp$ is supercuspidal: implicit in the theorem is the statement that $\Sigma|_{D_\wp}$ is independent of various choices made, as for instance, Harris and Taylor prove that the association $\Pi_\wp \to \Sigma|_{D_\wp}$ is purely local. This proves Theorem 1.

## 3  Acknowledgements

I would like to express my sincere thanks to Dipendra Prasad for his encouragement and his detailed comments on early versions of the manuscript: his help has been invaluable.

I would also like to thank Marie-France Vigneras for the interest she has shown in the present work, and for making available the preprint [V].

Finally I would like to thank National Univeristy of Singapore for support during May 1999 in which most of this paper got written.